\documentclass[11pt]{amsart}

\usepackage{amssymb}
\usepackage{amsmath}
\usepackage[small,nohug
]{diagrams}
\usepackage{amsthm, times, fullpage, float}
\usepackage{caption}
\usepackage{pstricks}
\usepackage{color}

\usepackage{palatino}

\newcommand{\Lg}{\mbox{$\mathfrak g$}}
\newcommand{\Ll}{\mbox{$\mathfrak l$}}
\newcommand{\Lh}{\mbox{$\mathfrak h$}}
\newcommand{\Lk}{\mbox{$\mathfrak k$}}
\newcommand{\Lp}{\mbox{$\mathfrak p$}}
\newcommand{\La}{\mbox{$\mathfrak a$}}

\newcommand{\Lm}{\mbox{$\mathfrak m$}}
\newcommand{\Ln}{\mbox{$\mathfrak n$}}

\newcommand{\Lt}{\mbox{$\mathfrak t$}}
\newcommand{\Lq}{\mbox{$\mathfrak q$}}
\newcommand{\Ls}{\mbox{$\mathfrak s$}}
\newcommand{\Lz}{\mbox{$\mathfrak z$}}

\newcommand{\Pf}{{\em Proof}. }
\newcommand{\EPf}{\hfill$\square$}
\newcommand{\Z}{\mbox{$\mathbb Z$}}
\newcommand{\R}{\mbox{$\mathbb R$}}
\newcommand{\C}{\mbox{$\mathbb C$}}

\newcommand{\SU}[1]{\mbox{$\mathrm{SU}(#1)$}}
\newcommand{\U}[1]{\mbox{$\mathrm{U}(#1)$}}
\newcommand{\SP}[1]{\mbox{$\mathrm{Sp}(#1)$}}
\newcommand{\SO}[1]{\mbox{$\mathrm{SO}(#1)$}}

\newcommand{\Spin}[1]{\mbox{$\mathrm{Spin}(#1)$}}
\newcommand{\F}{\mbox{$\mathrm{F}_4$}}
\newcommand{\G}{\mbox{$\mathrm{G}_2$}}

\newcommand{\su}[1]{\mbox{$\mathfrak{su}(#1)$}}
\newcommand{\Lu}[1]{\mbox{$\mathfrak{u}(#1)$}}
\newcommand{\ssp}[1]{\mbox{$\mathfrak{sp}(#1)$}}
\newcommand{\so}[1]{\mbox{$\mathfrak{so}(#1)$}}

\newtheorem{thm}{Theorem}

\newtheorem{cor}{Corollary}
\newtheorem{prop}{Proposition}
\newtheorem{lem}{Lemma}
\theoremstyle{remark}
\newtheorem{rem}{Remark}

\newtheorem*{notn*}{Notation}

\title{Tight Lagrangian homology spheres in\\
compact homogeneous K\"ahler manifolds}
\author{Claudio Gorodski and Fabio Podest\`a}

\address{Instituto de Matem\'atica e Estat\'\i stica, Universidade de
S\~ao Paulo, Rua do Mat\~ao, 1010, S\~ao Paulo, SP 05508-090, Brazil}
\email{gorodski@ime.usp.br}

\address{Dipartimento di Matematica e Informatica "Ulisse Dini", Universit\`a di Firenze, V.le Morgagni 67/A, 50100 Firenze, Italy}
\email{podesta@unifi.it}

\begin{document}

\begin{abstract}
For any irreducible compact homogeneous K\"ahler manifold,
we classify the compact tight Lagrangian submanifolds
which have the $\Z_2$-homology of a sphere.
\end{abstract}

\maketitle

\section{Introduction}

Let $M$ be a homogeneous K\"ahler manifold. Following~\cite{Oh},
we call a compact Lagrangian submanifold $L$ of $M$ \emph{globally tight} (resp.~\emph{locally tight} or simply~\emph{tight})
if the cardinality of the set $L\cap g\cdot L$ is equal to the sum of
$\Z_2$-Betti numbers of $L$, for every isometry $g$ of $M$
(resp.~every isometry
sufficiently close to the identity) such that the intersection is transversal.

It turns out that tightness has a bearing on the problem
of Hamiltonian volume minimization. For instance in~\cite{Oh} it is proved
that a tight Lagrangian submanifold of $\C P^n$ must be the
totally real embedding of $\R P^n$; an argument of Kleiner and Oh shows that
the standard $\R P^n$ in $\C P^n$ has the least volume among its Hamiltonian
deformations; Iriyeh~\cite{Ir} then notes that this gives uniqueness
of the Hamiltonian volume minimization problem for Hamiltonian deformations
of $\R P^n\subset \C P^n$ (similar results have been obtained
for the product of equatorial circles $S^1\times S^1\subset S^2\times S^2
=Q_2$~\cite{IOS,IS}).
More generally, real forms of Hermitian symmetric spaces
have recently been proved to be globally tight~\cite{TT}. The question
of classification of tight Lagrangian submanifolds in Hermitian
symmetric spaces was already posed in~\cite{Oh} and remains open.

Herein we take a different standpoint in that we allow $M$ to be
an arbitrary  compact homogeneous irreducible K\"ahler manifold
but we considerably restrict the topology of~$L$.
A compact homogeneous K\"ahler manifold $M$ is a
K\"ahler manifold
on which a compact connected Lie group of
isometries acts transitively. A simply-connected
compact homogeneous K\"ahler manifold $M$ is also called
a \emph{K\"ahlerian C-space}. In this case, it is known
that~$M$ is a homogeneous space $G/H$ where $G$ is a compact semi-simple Lie
group and $H$ is the centralizer of a toral subgroup of $G$
(in other words, it is a  \emph{(generalized) complex
flag manifold}); moreover $M$ is irreducible if and only if $G$ is a simple Lie group.
Our main result is:

\begin{thm}\label{main}
Let $M=G/H$ be a simply-connected irreducible compact homogeneous
K\"ahler manifold.
Let $L$ be a compact tight Lagrangian submanifold of~$M$.
Assume that $L$ has the $\Z_2$-homology of a sphere. \par
Then $L$ is an orbit of a compact subgroup of $G$,
and $M$ and $L$ are given, respectively up to  biholomorphic
homothety and up to congruence, as follows:
 \begin{itemize}
\item[(a)] $M$ is a complex quadric
$Q_n=\SO{n+2}/\SO2\times\SO{n}$ ($n\geq3$) and
$L\cong S^n$ is its standard real form, orbit of a subgroup
isomorphic to $\SO{n+1}$;
\item[(b)] $M$ is the twistor space
$Z=\SU{n+1}/\mathrm{S}(\U 1\times \U 1\times \U {n-1})$ ($n\geq 3$)
of the complex Grassmannian of $2$-planes
$\mbox{Gr}_2(\C^{n+1})$ endowed with its standard K\"ahler-Einstein
structure and $L\cong S^{2n-1}$ is an orbit of a subgroup
isomorphic to $\U{n}$;
\item[(c)] $M$ is the full flag manifold
$\SU 3/{\mathrm T}^2$ endowed with its K\"ahler-Einstein
homogeneous metric and $L\cong S^3$ is an orbit of a subgroup
isomorphic to~$\U2$;
\item[(d)] $M = \SP{n+2}/\U 2\times \SP n$ ($n\geq 1$) and
$L\cong S^{4n+3}$ is an orbit of a subgroup isomorphic to
$\SP 1 \times \SP{n+1}$;
\item[(e)] $M = \F/{\mathrm T}^1\cdot\Spin 7$ and
$L\cong S^{15}$ is an orbit of a subgroup isomorphic to $\Spin 9$.
\end{itemize}
Moreover $L$ coincides with a connected component of the fixed point set of an antiholomorphic isometric involution of $M$.
\end{thm}

\begin{rem} Note that the space $\SU 3/{\mathrm T}^2$ admits only one invariant complex structure up to equivalence, while the spaces appearing in (d) and (e) carry only one invariant K\"ahler structure up to biholomorphism and homothety, so that it is not necessary to specify which structure we are considering. For further details we refer to \S 2.
\end{rem}

In section~\ref{oh},
we briefly review some basic facts about
homogeneous K\"ahler manifolds. The proof
of Theorem~\ref{main} is scattered throughout sections~\ref{lag},
\ref{og} and~\ref{final}. In particular in \S~\ref{final} we also show that a real flag manifold can be always embedded as a tight real form of a suitable complexification given by a complex flag manifold, see Proposition \ref{real flags}.

\begin{notn*}
For a compact Lie group, we denote its Lie algebra by the corresponding lowercase gothic letter. If a group $G$ acts on a manifold $M$, for every $X\in \Lg$ we denote by $X^*$ the corresponding vector field on $M$ induced by the
$G$-action.
\end{notn*}

\section{Preliminary material}\label{oh}

Let
$M=G/H$ be a generalized flag manifold, where
$G$ is a compact connected semisimple Lie group and $H$ is the centralizer
of a toral subgroup of $G$. We shall recall the standard description of
invariant K\"ahler structures on $M$ (see e.g.~\cite{BFR,Alek}).

Denote by $p$ the basepoint and by $\langle,\rangle$ the
negative of the Cartan-Killing form of~$\Lg$.
Then there is a reductive decomposition $\Lg = \Lh \oplus \Lm$ where $\Lm$ is the
orthogonal complement of $\Lh$, and we can as usual identify
$\Lm\cong T_{p}(G/H)$ via $X\mapsto X^*_p$. Since $\Lh$ is of maximal rank in $\Lg$,
there is a maximal Abelian subalgebra $\Ls$ of $\Lg$ contained
in~$\Lh$. Then the complexification $\Ls^{\mathbb C}$ is a Cartan subalgebra of
$\Lg^{\mathbb C}$ and we denote by $\Delta$ the corresponding root system.
Each root space of $\Lg^{\mathbb C}$ is either contained in
$\Lh^{\mathbb C}$ or in $\Lm^{\mathbb C}$ and thus
there is an associated partition $\Delta=\Delta_H\cup\Delta_M$.
Define the real subspace
$$ \Lt = i\ \Lz(\Lh) \subset \Ls^{\mathbb C}$$
where $\Lz(\Lh)$ is the center of $\Lh$. Every root $\alpha\in \Delta$ is real valued when restricted to $\Lt$ and
the restriction $\alpha|_{\Lt}\in \Lt^*$ is called a \emph{$T$-root} (note that the elements of $\Delta_H$ restrict to zero).
Note that the set $\Delta_T\subset \Lt^*$ of all $T$-roots is not a root system.
Its significance for us lies in the fact that there is a natural
bijective correspondence between the set of $G$-invariant complex structures
on $M$  and the set of $T$-chambers in $\Lt$, where
a \emph{$T$-chamber} is a connected component of the regular set $\Lt_{reg}$
of $\Lt$, namely, the complement of the union
of the hyperplanes $\ker\lambda$ for~$\lambda\in\Delta_T$. For later
reference, we recall that in case $\dim\Lt=1$ or
$\Lh$ is Abelian, any two $G$-invariant
complex structures on $M$ are biholomorphic~\cite[13.8]{BH} (see also~\cite[p.57]{Nis}).

On the other hand, there is a natural bijective correspondence between
the set of $G$-invariant symplectic structures $\omega$ on $M$ and the set of
regular elements $\xi\in i\ \Lt_{reg}$ given by
\[ \omega_p(X^*_p,Y^*_p) = \langle[X,Y],\xi\rangle \]
where $X$, $Y\in\Lg$. Finally the $G$-invariant K\"ahler metrics on $M$ are all
given by
\[ g_p=\omega_p(\cdot,J\cdot) \]
where $\omega$ is the invariant symplectic
structure associated to $\xi\in i\ \Lt_{reg}$ and
$J$ is the invariant complex structure associated to the $T$-chamber
containing~$-i\ \xi$.

Henceforth we fix an invariant K\"ahler structure on $M$.
Since $\Lh$ coincides with the centralizer of~$\xi$ in~$\Lg$,
there is a canonical embedding of $M$ into $\Lg$ as the adjoint orbit
$\mathrm{Ad}(G)\cdot\xi$  mapping $p$ to $\xi$.
Now for each~$X\in\Lg$, the vector field $X^*$ on
$M$ is given by $X^*_q=[X,q]$, where $q\in M$, and it is easily seen that
\[ \omega_q(X^*_q,Y^*_q):=\langle [X,Y],q\rangle, \]
showing that $\omega$ coincides with the Kirillov-Kostant-Souriau
symplectic structure. It follows that~$X^*$ is
Hamiltonian with corresponding
potential function given by
the height function $h_X(q)=\langle q,X\rangle$ for $q\in M$.
Thus the moment map
\[ \mu : M \to\Lg^*\cong\Lg,\quad \mu(q)(X)=h_X(q) \]
is just the inclusion (above we have identified $\Lg$ with its dual via the
Cartan-Killing form).

For later reference, we quote
the following result due to Onishchik~\cite[p.~244]{On}.

\begin{thm}\label{oni} If $G$ is a compact connected simple Lie group
and acts effectively on $M=G/H$, then $G$ coincides,
up to covering, with the identity
component of the full isometry group $Q$,
with the following exceptions:
\begin{itemize}
\item[(a)] $M=\C P^{2n+1}$ and $\Lg = \ssp {n+1}$, $\Lh = \Lu{1} \oplus \ssp n$, $\Lq = \su {2n+2}$ ;
\item[(b)] $\Lg = \so{2n-1}$, $\Lh = \Lu {n-1}$, $\Lq = \so {2n}$, \ $n\geq 4$;
\item[(c)] $M = Q_5$ and $\Lg = \Lg_2$, $\Lh = \Lu 2$, $\Lq = \so 7$.
\end{itemize}
\end{thm}

\section{Lagrangian submanifolds}\label{lag}

We keep the notation from the previous section and
consider a compact Lagrangian submanifold $L$ of $M$ through the basepoint~$p$.
Denote by $K$ the identity component of the stabilizer subgroup of $L$ in $G$.
Then $K$ is a closed subgroup of $G$ that acts effectively on $L$ by the Lagrangian
property of $L$.
Define the linear map
\begin{equation}\label{norm}
 \sigma: \Lg\to\Gamma(\nu L),\quad X\mapsto(X^*|_L)^\perp,
\end{equation}
where $\Gamma(\nu L)$ is the space of sections
of the normal bundle $\nu L$ and $()^\perp$
denotes the normal component to $L$.
The map~(\ref{norm}) is clearly $K$-equivariant;
note that its kernel coincides with the Lie algebra $\Lk$ of $K$.
Choose a reductive complement $\Lk^\perp$ and
$K$-equivariantly identify
$\Lk^\perp$ with the image~$V$ of~$\sigma$. Now we can write
\begin{equation}\label{1}
 \Lg = \Lk + V.
\end{equation}

\begin{lem}\label{2}
For $X\in\Lg$, the critical points of the height function
$h_X|_L$ are precisely the zeros of the vector field
$\sigma(X)\in V$.
\end{lem}

\Pf A point $q\in L$ is a critical point of $h_X|_L$ if and only
if $X$ is perpendicular to $T_qL$. Given $v\in T_qL$,
there exists $Y\in\Lg$ such that $Y^*_q=v$
and then $\omega_q(X^*_q,v)=\langle [X,Y],q\rangle=\langle X,[Y,q]\rangle
=\langle X,v\rangle$.
Now $q$ is a critical point of $h_X|_L$ if and only if
$\omega_q(X^*_q,T_qL)=0$. Since $L$ is Lagrangian,
this implies that $X^*_q\in T_qL$, as desired. \EPf

\medskip

Let $\varphi:L\to V$ be the restriction of the orthogonal projection
$\Lg\to V$. We elaborate on an idea of Oh, use the previous lemma to translate
the tight Lagrangian property of~$L$ to the taut property of~$\varphi(L)$, and apply
known results to the latter.

\begin{prop}\label{taut}
The map $\varphi:L\to V$ is a $K$-equivariant full embedding.
Moreover, if $L$ is tight Lagrangian in $M$ then
$\varphi(L)$ is a taut submanifold of the Euclidean space $V$;
the converse holds in the case~$G$ coincides,
up to covering, with the
identity component of the isometry group of~$M$.
\end{prop}

\Pf Since~(\ref{1}) is a reductive decomposition,
it is clear that $\varphi$ is $K$-equivariant.
The moment map of the restricted $K$-action on $M$ is
$\pi_{\mathfrak k}\circ\mu$, where $\pi_{\mathfrak k}:\Lg\to\Lk$
is the orthogonal projection; since~$L$ is a $K$-invariant
Lagrangian submanifold, we have
\[ \pi_{\mathfrak k}\circ\mu(L) = \eta, \]
where $\eta$ is a constant central element of $\Lk$.
Denote by $\pi_V:\Lg\to V$ the orthogonal projection. Then
\begin{equation}\label{image} \mu|_L=(\pi_{\mathfrak k}+\pi_V)\circ\mu|_L = \eta + \varphi. \end{equation}
Since $\mu|_L$ is the inclusion into~$\Lg$, this shows that
$\varphi$ is an embedding.

If $\varphi:L\to V$ is not full, then
$\varphi(L)$ is contained in an affine hyperplane of $V$, namely,
$\langle\varphi(q),\zeta\rangle=\langle q,\zeta\rangle$ is a constant for every $q\in L$
and some nonzero $\zeta\in V$. This implies that the height function
$h_\zeta$ is constant on $L$ and thus, by Lemma~\ref{2},
$\sigma(\zeta)$ is
 everywhere zero, namely~$\zeta\in \Lk$, a
contradiction to~$\Lk\cap V = \{0\}$. This proves that $\varphi:L\to V$ is full.

Recall that $\varphi:L\to V$ is by definition a tight embedding
if and only if every height function $h_X|_{\varphi(L)}$
for $X\in V$ which is a Morse function
is also perfect, i.e.~has the minimum number of critical
points allowed by the Morse inequalities~\cite{C-C}. By Lemma~\ref{2},
this is equivalent to $X^*|_{L}$ having a number of zeros
equal to the sum of $\Z_2$-Betti numbers of $L$ for generic
$X\in V$. Since such $X$ are the infinitesimal generators
of one-parameter groups of isometries of $M$, the latter condition
follows from the tightness of $L$ in $M$,
and is equivalent to it in case $\Lg$ coincides with the Lie algebra of all
isometries of $M$.

Finally, note that $M$ is a $G$-orbit so it is contained
in a round sphere of $\Lg$.
By~(\ref{image}) also $\varphi(L)$ is contained in a round sphere of $V$.
A submanifold contained in a round sphere in a Euclidean space is tight if and only if all
distance functions are perfect Morse functions, namely,
if and only if it is taut~\cite{C-C}.
Furthermore, in this situation the set of critical
points of a distance function will also occur as
the set of critical points of a height function, and vice versa. \EPf

\begin{rem}\label{zk=0}
It follows from~(\ref{image}) that $L\subset V$ if $\Lk$ is centerless.
\end{rem}

\begin{rem}\label{hsp}
If $K$ is a symmetric subgroup of $G$, then its orbits in $V$ are taut
submanifolds (see e.g.~\cite{GTh2}). Thus it easily follows
from Proposition~\ref{taut} that real forms of Hermitian symmetric
spaces of compact type are locally tight; we omit the details.
Note that this result already follows from the work of Takeuchi and
and Kobayashi~\cite{TK}. Moreover it has been recently proved
by Tanaka and Tasaki that those real forms are indeed globally tight~\cite{TT}.
\end{rem}

Recall that the Chern-Lashof theorem~\cite{C-L}
implies that a taut and substantial smooth embedding of a
$\Z_2$-homology sphere
into an Euclidean sphere must be round and have codimension one (see also~\cite{No-Ro}).
Hence:

\begin{cor}\label{dimV1}
If $L$ is a compact tight Lagrangian $\Z_2$-homology sphere then
$\varphi(L)$ is a codimension one round sphere in $V$.
In particular
\[  \dim V = \dim L +1. \]
\end{cor}

\section{The Ohnita-Gotoh formula}\label{og}

Keep the notation from the previous
two sections and assume for the moment that
$L$ is an arbitrary compact Lagrangian submanifold of~$M$.
We introduce the subspace $\Ll$ of $\Lm$ corresponding to~$T_pL$.
The $G$-invariant Riemannian metric on $M$ corresponds to
an $\mathrm{Ad}(H)$-invariant
inner product
in $\Lm$; let $\Ll^\perp$ be the orthogonal complement of $\Ll$
in $\Lm$. Also, denote the normalizer subalgebra of~$\Ll$ in~$\Lh$ by~$\Ln$.

The following proposition elaborates on results by Ohnita~\cite{Ohn} and Gotoh~\cite{Go}.

\begin{prop}\label{dimV2} We have
\begin{equation}\label{ineq}\dim V \geq \dim L + \dim \Lh - \dim \Ln.\end{equation}
Moreover, if equality holds then $L$ is homogeneous under the action of $K$ and $\Ln\subset \Lk \cap \Lh$.
\end{prop}

\Pf Throughout we identify $T_pL\cong \Ll$ and
$\nu_pL\cong \Ll^\perp$ whenever clear from context.
 We consider the diagram
\begin{equation}\label{diag}
\begin{diagram}
\Lg &\rTo^{\scriptstyle\Psi=\Psi_1\oplus\Psi_2}&
\Ll^\perp\oplus\mathrm{Hom}(\Ll,\Ll^\perp) \\
\dTo_{\scriptstyle\sigma} && \dTo_{\scriptstyle\cong}\\
V &\rTo_{\scriptstyle\Phi_{p}}& \nu_{p}L\oplus\mathrm{Hom}(T_{p}L,\nu_{p}L)
\end{diagram}
\end{equation}
where $ \Psi_1:\Lg\to\Ll^\perp $ is the projection with respect to
the vector space direct sum decomposition $\Lg = \Lh + \Ll + \Ll^\perp$, the map
$\Psi_2:\Lg\to\mathrm{Hom}(\Ll,\Ll^\perp)$ is given by
\[ \Psi_2(X)(Y) = \left(\nabla_{Y^*}X_\mathfrak{m}^*\right)|_p^\perp + [X_\mathfrak{h},Y]^\perp - B(X_{\mathfrak{l}},Y), \]
where $B:\Ll\times\Ll\to\Ll^\perp$ is the second fundamental form of $L$ in $M$ at~$p$,
and
\[ \Phi_{p}(\eta)=(\eta_{p},\nabla^\perp\eta|_{p}) \]
for $X\in \Lg$, $Y\in \Ll$ and $\eta=(X^*|_L)^\perp\in V$.

The commutativity of diagram \eqref{diag} follows from $\sigma(X)|_p=(X^*_p)^\perp$ and
\begin{equation}\label{comm}\left(\nabla_{Y^*}X_\mathfrak{m}^*\right)^\perp|_p +
[X_\mathfrak{h},Y]^\perp - B(X_{\mathfrak{l}},Y) = \nabla^\perp_{Y^*}\sigma(X)|_{p}
\end{equation}
for $X\in \Lg$ and $Y\in \Ll$. In turn, we check~(\ref{comm}) as follows.
\begin{eqnarray*}
\left(\nabla_{Y^*}X_\mathfrak{m}^*\right)^\perp|_p &=&
\left(\nabla_{Y^*}X^*\right)^\perp|_p -\left(\nabla_{Y^*}X_\mathfrak{h}^*\right)^\perp|_p \\
&=&\left(\nabla_{Y^*}(X^*)^\top\right)^\perp|_p
+\left(\nabla_{Y^*}(X^*)^\perp\right)^\perp|_p
-\left(\nabla_{Y^*}X_\mathfrak{h}^*\right)^\perp|_p \\
&=&B(X_{\mathfrak l},Y)+\nabla_{Y^*}^\perp\sigma(X)|_p-(\nabla_{Y^*}X_{\mathfrak h}^*)^\perp|_p.
\end{eqnarray*}
Finally, the result follows from the formula
\[ \left(\nabla_{W^*}U^*\right)|_p = [U,W]^*_p \]
for $U\in \Lh$ and $W\in \Lm$, which is easily proved using~(7.27)
in~\cite{Be} and the fact that
$\mathrm{ad}(U)\in \mathrm{End}(\Lm)$ is skew-symmetric
with respect to the metric in $\Lm$.

It follows from the commutativity of diagram~\eqref{diag} and the surjectivity of $\sigma$ that $\mathrm{im}(\Psi) = \mathrm{im}(\Phi_p)$ and thus
\begin{equation}\label{ineq1}
\dim V \geq \dim \mathrm{im}(\Phi_p) = \dim \mathrm{im}(\Psi).
\end{equation}
It is obvious that $\Psi(\Ll^\perp)\cap \Psi(\Lh + \Ll) = \{0\}$.
Therefore
\begin{eqnarray}\label{ineq2}
\dim \mathrm{im}(\Psi) &=& \dim \Psi(\Ll^\perp) + \dim \Psi(\Lh + \Ll) \\ \nonumber
&\geq& \dim \Ll^\perp + \dim \Psi(\Lh)\\ \nonumber
&=& \dim L + \dim \Lh - \dim \ker(\Psi|_{\Lh})\\ \nonumber
& =& \dim L + \dim \Lh - \dim \Ln,
\end{eqnarray}
proving~\eqref{ineq}.

In the case of equality in~(\ref{ineq}), we follow~\cite{Go}. Use
(\ref{ineq1}) and~(\ref{ineq2}) we see that $\Phi_p$ is injective and
$\Psi(\Lh + \Ll) = \Psi(\Lh)$. Now for given~$v\in T_pL$
we can find $X\in \Ll$ and $Y\in \Lh$
with $X^*_p=v$ and $\Psi(Y) = -\Psi(X)$.
Therefore $\sigma(X+Y)=0$, namely, $X+Y\in \Lk$.
This proves that $K$ acts transitively on $L$.
Moreover $\Ln = \ker(\Psi|_{\Lh}) \subset \ker\sigma = \Lk$. \EPf

\begin{cor}\label{L-homogeneity}
If $L$ is a compact tight Lagrangian $\Z_2$-homology sphere
then $L$ is homogeneous under $K$ and
$\Ln=\Lk\cap\Lh$ is the isotropy subalgebra of~$\Lk$ at~$p$
and a codimension one ideal of~$\Lh$.
Moreover $(\Lg,\Lk)$ is either a symmetric pair of rank one
or~$(\mathfrak g_2, \su 3)$ or~$(\so 7, \mathfrak g_2)$.
\end{cor}

\Pf We see that $\Ln\subsetneq\Lh$ by noting that
$\xi\not\in\Ln$. Indeed if $[\xi,\Ll]\subseteq \Ll$ then
\begin{equation}\label{zero} 0 = \omega_p(T_pL,T_pL) = \langle [\xi,\Ll],\Ll\rangle, \end{equation}
which implies $[\xi,\Ll] = \{0\}$, contradicting the facts that
the centralizer of $\xi$ in $\Lg$ is $\Lh$, and $\Lh\cap\Ll = \{0\}$.

Further, it follows from~Corollary~\ref{dimV1} and
Proposition~\ref{dimV2} that
$\dim\Ln=\dim\Lh-1$, $L$ is $K$-homogeneous and $\Ln\subset\Lk\cap\Lh=\Lk_p$.
The reverse inclusion $\Lk_p\subset\Ln$ is obvious and therefore $\Ln=\Lk\cap\Lh$.

It also follows from Corollary~\ref{dimV1} that
$K$ acts on $V$ with cohomogeneity one and the last claim
follows (see e.g.~\cite[3.12]{HPTT}).
\EPf

\section{End of the proof of Theorem~\ref{main}}\label{final}

First we explain a standard construction which allows one to
construct compact tight Lagrangian submanifolds in
suitable complex flag manifolds.

Let $K$ be a connected \emph{symmetric} subgroup of a compact
connected semisimple Lie group $G$ and consider the
decomposition $\Lg = \Lk + V$ into eigenspaces of the involution.
Any orbit of $K$
on $V$, say $L=\mathrm{Ad}(K)\cdot\xi$ for some $\xi\in V$,
is called a \emph{(generalized) real flag manifold}.
There is a natural ``complexification'' of $L$,
namely, we next show that the adjoint orbit $M=\mathrm{Ad}(G)\cdot\xi$
is a complex flag manifold containing $L$
as the connected component of the fixed point
set of an anti-holomorphic involutive isometry.

Since $G$ is compact,
it embeds into its complexification $G^{\mathbb C}$ as a maximal
compact subgroup. The Lie algebra of $G^{\mathbb C}$ is
$\Lg^{\mathbb C}=\mathbb \Lg\otimes\mathbb C$ and admits
$\Lg_0 = \Lk + iV$ as a non-compact real form;
let $G_0$ denote the corresponding connected subgroup
of $G^{\mathbb C}$ and $\tau$ the associated conjugation of $G^{\mathbb C}$
over~$G_0$.

Fix a maximal Abelian subalgebra $\La$ of $iV$ containing $a:=i\xi$
and consider the restricted root decomposition
$\Lg_0=Z_{\mathfrak k}(\La)+\La+\sum_{\lambda\in\Sigma}\Lg_{0,\lambda}$
where $Z_{\mathfrak k}(\La)$ denotes the centralizer of $\La$
in $\Lk$. Choose a positive restricted root system
$\Sigma^+\subset\Sigma$ so that
$\Lg_0 = \Lk+\La+\Ln$ is an Iwasawa decomposition,
where $\Ln=\sum_{\lambda\in\Sigma^+}\Lg_{0,\lambda}$.
As a homogeneous space, $L=K/Z_K(a)$, where $Z_K(a)$
is the centralizer of $a$ in $K$ and its Lie algebra
is $Z_{\mathfrak k}(a)=Z_{\mathfrak k}(\La)+
\sum_{\lambda\in\Sigma^+\atop\lambda(a)=0}\Lk_{\lambda}$,
where $\Lk_{\lambda}=(\Lg_{0,\lambda}+\Lg_{0,-\lambda})\cap\Lk$.

It turns out that $G_0$ acts on $L$. To see that,
recall that a minimal parabolic subalgebra of $\Lg_0$ is
any subalgebra conjugated to $\Lp_{0,min}=Z_{\mathfrak k}(\La)+\La+\Ln$,
and a parabolic subalgebra of $\Lg_0$ is any subalgebra
containing a minimal parabolic subalgebra (see e.g.~\cite[\S 1.2.3 and~1.2.4]{warner}). Now
\[ \Lp_0:=Z_{\mathfrak k}(\La)+\La+\Ln+
\sum_{\lambda\in\Sigma^+\atop\lambda(a)=0}\Lg_{0,-\lambda} \]
is a parabolic subalgebra of $\Lg_0$.
The normalizer $P_0$ of $\Lp_0$ in $G_0$ is called
a parabolic subgroup.
Let $\Theta$ be the set of simple restricted roots
$\lambda$ satisfying $\lambda(a)=0$, and let
$\La_{\Theta}$ be
the subspace of $\La$ that $\Theta$ annihilates.
By Theorem~1.2.4.8 in~\cite{warner},
$P_0=M_{\Theta}AN$ where $A=\exp(\La)$, $N=\exp(\Ln)$ and
$M_\Theta$ is the centralizer $Z_K(\La_{\Theta})$ of $\La_{\Theta}$ in $K$
(loc.~cit., p.~73).
Note that $a$ is a generic element in $\La_{\Theta}$, so
$Z_K(\La_{\Theta})=Z_K(a)$. In particular,
$K\cap P_0=M_\Theta=Z_K(a)$. The group $K$ acts
by left translations on $G_0/P_0$ with an orbit that is open
(by counting dimensions, since $\dim\Lk_{\lambda}=\dim\Lg_{\lambda}$)
and closed (by compactness of $K$), so
$K/Z_K(\xi)=K/Z_K(a)=K/K\cap P_0=G_0/P_0$.
This realizes the real flag $L$ as a $G_0$-homogeneous space.

On the other hand, $G^{\mathbf C}$ acts on $M$. Indeed, let $\Lt_{\mathfrak k}$
be a maximal Abelian subalgebra of $Z_{\mathfrak k}(\La)$. Then
$\Ls=\Lt_{\mathfrak k}+\La$ is a Cartan subalgebra of $\Lg_0$,
and $\Ls^{\mathbb C}$ is a Cartan subalgebra of $\Lg^{\mathbb C}$
with root system $\Delta$ and root decomposition
\[ \Lg^{\mathbb C}=\Ls^{\mathbb C}+ \sum_{\alpha\in\Delta}\Lg_\alpha^{\mathbb C}. \]
The roots are real valued on $\Ls_{\mathbb R}=i\Lt_{\mathfrak k}+\La$, and we
take a lexicographic order that takes $\La$ before
$i\Lt_{\mathfrak k}$. The point is that a restricted root of the form
$\lambda=\alpha|\La$ for $\alpha\in\Delta$ is positive if and only if
$\alpha\in\Delta^+$. Now
\begin{eqnarray*}
\Lp &:=& \Lp_0\otimes\mathbb C \\
    &=& Z_{\mathfrak k}(\La)^{\mathbb C} + \La^{\mathbb C} +\Ln^{\mathbb C} +
\sum_{\lambda\in\Sigma^+\atop\lambda(a)=0}\Lg_{0,-\lambda}\otimes\mathbb C \\
&=&\left(\Lt_{\mathfrak k}^{\mathbb C}+\sum_{\alpha|\mathfrak a=0}\Lg_{\alpha}^{\mathbb C}\right)
+\La^{\mathbb C}+\sum_{\alpha\in\Delta^+\atop\alpha|\mathfrak a\neq0}\Lg_{\alpha}^{\mathbb C}
+\sum_{\alpha\in\Delta^+\atop\alpha|\mathfrak a\neq0, \alpha(a)=0}\Lg_{-\alpha}^{\mathbb C}\\
&=&\underbrace{\Ls^{\mathbb C} + \sum_{\alpha\in\Delta^+}\Lg_{\alpha}^{\mathbb C}}_{\textrm{Borel
subalgebra}}+
\sum_{\alpha\in\Delta^+\atop\alpha(a)=0}\Lg_{-\alpha}^{\mathbb C}
\end{eqnarray*}
is a parabolic subalgebra of $\Lg^{\mathbb C}$, since it contains
a Borel subalgebra of $\Lg^{\mathbb C}$. The normalizer~$P$
of $\Lp$ in $G^{\mathbb C}$ is called a parabolic subgroup of $G$.
It is closed and connected. As in the real case above,
it follows (even easier) that $G\cap P=Z_G(\xi)$
(cf.~\cite[Corollary~2.7]{wolf})
and $G/Z_G(\xi)=G/G\cap P=G^{\mathbb C}/P$. This realizes the complex flag
$M$ as a $G^{\mathbb C}$-homogeneous space.

The involution $\tau$ of $G^{\mathbb C}$ stabilizes $P$ so
induces an involution $\bar\tau$ of $M$
whose connected component through the basepoint
coincides with~$L$~\cite[Theorem~3.6]{wolf}.

Consider the Kirillov-Kostant-Souriau invariant symplectic
form on $M$ defined at $\xi$ by
\[ \omega_\xi(X,Y)=\langle [X,Y],\xi\rangle \]
for $X$, $Y\in T_\xi M$. Since $\tau$ preserves the Killing form of $\Lg$
and maps $\xi$ to $-\xi$, we see that $\bar\tau$ is antisymplectic
and thus $L$ is a Lagrangian submanifold of $M$. Note that $\bar\tau$
is also antiholomorphic with respect to the complex structure $J$
on $M$ given by the $T$-chamber containing $\xi$, and hence
it is an isometry with respect to the K\"ahler metric~$g$
induced by $\omega$ and $J$.

The real flag manifold $L$ is called a \emph{real form} of the complex flag
manifold~$M$ endowed with the K\"ahler metric $g$. In case $L$ is already
a complex flag manifold viewed
as a real flag manifold (namely, $G_0$ is a complex semisimple Lie group
viewed as a real Lie group), $M$ can be identified
with $L\times \bar L$, where $\bar L$ is equipped with the opposite
complex structure and $L$ sits in $M$ as the diagonal.

\begin{rem}
Note that we can start with any real flag manifold and ``complexify''
it to a complex flag manifold. Conversely, if we start
with a complex flag $G^{\mathbb C}/Q$ and fix a real form $G_0$, it is not always
true that there is a $G_0$-orbit in $G^{\mathbb C}/Q$ which is a real form;
it is not true, for instance, for a full flag $G^{\mathbb C}/B$
where $G^{\mathbb C}$ is a complex semisimple Lie group
with Lie algebra $\Lg^{\mathbb C}$, and $B$ is a Borel subgroup
such that its Lie algebra contains the complexification of
a maximally split Cartan subalgebra $\Lt_{\mathfrak k}+\La$
of a real form $\Lg_0$, and $\La$ contains no regular element
of $\Lg^{\mathbb C}$
\cite[p.~1139]{wolf}. On the other hand,
if $G_0$ is the Cartan normal real form,
we can always find a $G_0$-orbit in $G^{\mathbb C}/Q$ which is a real form.
\end{rem}

As far as we know, real forms of complex flags manifolds
other than Hermitian symmetric spaces have not been
explicitly classified (see~\cite[Theorem~3.6]{wolf}, though).
Compare the next result with Remark~\ref{hsp}.

\begin{prop}\label{real flags}
Consider the real form  $L$ of the complex flag manifold $M$ constructed
above. Then $L$ is a tight Lagrangian submanifold of $M$.
\end{prop}

\Pf We need only prove the tightness. A symmetric space of compact type
splits into the direct product of irreducible symmetric spaces
of compact type, and the linear isotropy representation splits
accordingly, so we may assume $G/K$ is irreducible.

Suppose first $G/K$ is type of I, namely, $G$ is simple.
Since~$K$ is a symmetric subgroup of $G$, its orbits in~$V$
are taut submanifolds~\cite{GTh2}.
If the pair $(\Lg,\Lh)$ is not listed in Theorem~\ref{oni},
Proposition~\ref{taut} already implies that $L$ is tight.
Otherwise, the center of $\Lh$ is one-dimensional and therefore
$M = G/H$ admits precisely one invariant complex structure up to
biholomorphism and one invariant compatible K\"ahler metric up to homothety.
The complex manifold $M = G/H = G^\mathbb C/P$ can also be written as
$M = Q/H'$, where $Q$ denotes the connected group of all biholomophisms of
$M$ and the space $Q/H'$ is a Hermitian symmetric space.
The symmetric metric $\tilde g$ on $Q/H'$ (unique up to homothety)
is also $G$-invariant and is a scalar multiple of $g$. This means that the
submanifold $L$ is a real form of $M$ also with respect to
the symmetric metric $\tilde g$, and therefore it is tight
(see Remark~\ref{hsp}).

In case $G/K$ is of type II, $M=L\times\bar L$ (see above)
and the proof is analogous. \EPf

\medskip

The case of interest for us is that in which $(\Lg,\Lk)$ has rank one.
Here $L$ is a sphere, umbilic in $V$, for $\xi\neq0$.

Conversely, we now proceed to the proof of Theorem~\ref{main}
and examine to which extent the above construction supplies the
possible examples.
So let $M=G/H$ and $L$ be as in the statement of the theorem
and write $\Lg=\Lk+V$ as in section~\ref{lag}. We view
$M$ as the adjoint orbit through $\xi\in\Lg$ and assume $\xi\in L$.
Corollary~\ref{L-homogeneity} says that $L$ is homogeneous
under $K$, $\Lh=\Lk_\xi\oplus\mathfrak{u}(1)$
and $(\Lg,\Lk)$ is a symmetric pair of rank one or
one of two other pairs.

Note that if $\Lk$ is centerless, then we must have
$\xi\in V$ by Remark~\ref{zk=0}; if, in addition, the dimension
of the center of $\Lh$ is one, then $\xi$ lies in the
$\mathfrak{u}(1)$-summand of $\Lh$ and the K\"ahler structure on $M$ is unique,
up to biholomorphic homothety.
This is the situation for the symmetric pairs
$(\so {n+2},\so {n+1})$ ($n\geq 3$),
$(\ssp {n+2},\ssp{1}\oplus \ssp {n+1})$ ($n\geq1$) and
$(\mathfrak f_4,\so 9)$ (in these cases $(\Lg,\Lh)$
is not listed in Theorem~\ref{oni}), for which the construction above yields
the examples described in parts~(a), (d) and~(e) in the
statement of the theorem.

In case $(\Lg,\Lk)=(\Lg_2,\su3)$, the center of $\su3$ is zero and
that of $\Lh=\mathfrak{u}(2)$ is one-dimensional, so we
also must have $\xi\in V$. Since $\su3\subset\Lg_2$ is spanned
by the long roots of $\Lg_2$, $\Lh$ is the centralizer of a short
root and thus $M=\G/\U2$ is again the quadric $Q_5=\SO 7/\SO2\times \SO5$
and $L\cong\SU3/\SU2\cong S^5$.

In case $(\Lg,\Lk)=(\so7,\mathfrak g_2)$
again $\xi\in V$, so $\Lh\cong\Lu3$ is the centralizer of a
short root of $\so7$. It is known that $M= \SO 7/\U 3 \cong \SO 8/\U 4$ and an
outer automorphism $\tau$ of $\so 8$ induces a diffeomorphism between $M$ and,
again, the quadric $Q_6 = \SO 8/\SO 2 \times \SO 6$. Since $\G$ acts
with cohomogeneity one on $\SO7/\U3$~\cite[p.573]{K}, $\G$ and $\SO7$
share the same orbits in $Q_6$ and $L$ is congruent to the standard real form~$S^6$.

The last case we need to consider is $(\Lg,\Lk) = (\su {n+1},\mathfrak{s}(\mathfrak{u}(1) \oplus \mathfrak{u}(n)))$ ($n\geq2$). This case is somewhat more involved
because the center of $\Lk$ is non-trivial and $\Lh\cong \mathfrak{u}(1)\oplus \mathfrak{u}(1)\oplus \su{n-1}$ has
a center of dimension bigger than one.
We distinguish between two cases, namely, $n\geq 3$ and~$n = 2$.

\smallskip

(a) $n\geq 3$. The flag manifold $M$ is known to admit precisely two inequivalent $\SU{n+1}$-invariant complex structures (see e.g.~\cite{Nis}).
We note that
the tangent space $T_\xi M$ splits under the isotropy representation of $H=\mathrm{S}(\U 1\times\U 1\times\U {n-1})$ as the sum  $\bigoplus_{i=1}^3 V_i$ where
$V_i$ are mutually inequivalent complex $H$-submodules with $\dim_{\mathbb C} V_1 = 1$ and $\dim_{\mathbb C} V_2 = \dim_{\mathbb C} V_3 = n-1$.
It is relatively easy
to apply the machinery of section~\ref{oh} to show
that the two inequivalent invariant complex structures $J_0$, $J_1$
on $M$ can be described as follows (compare~\cite[\S13.9]{BH}):
if $(u,z,w)\in \bigoplus_{i=1}^3 V_i$, then
$J_\alpha (u,z,w) = (iu,iz,(-1)^\alpha iw)$ for $\alpha =0,1$.

We now consider the standard $\SU{n+1}$-equivariant fibration $\pi:M\to W$, where $W$ denotes the complex
Grassmannian $\mbox{Gr}_2(\C^{n+1}) = \SU{n+1}/\mathrm{S}(\U 2\times \U{n-1})$. It is clear that the projection $\pi$ is holomorphic when we endow $M$ with the complex structure $J_0$ and $W$ with its standard complex structure $J_W$ as an Hermitian symmetric space. On the other hand, the space $W$ is also a homogeneous quaternion-K\"ahler manifold, a so called Wolf space, endowed with an invariant quaternion K\"ahler structure $Q$, namely a rank three subbundle  $Q\subset \mbox{End}(TW)$ locally spanned by three local complex structures $\{I,J,K\}$ with $IJK = -\mbox{Id}$. It is well known that $J_W$ is not even a local section of $Q$ (see e.g.~\cite[14.53(b)]{Be}) and therefore $(M,J_1)$ is biholomorphic to the twistor space $Z$ of $W$.

\begin{lem}\label{twistor} If $M$ has an invariant K\"ahler structure $(g,J)$ admitting a compact tight Lagrangian $\mathbb Z_2$-homology sphere $L$, then $J$ is equivalent to $J_1$ and $g$ is K\"ahler-Einstein, i.e. $(M,g,J)$ is biholomorphically homothetic to the twistor space of $W$.
\end{lem}
\Pf We know there is a subgroup $K\subset \SU{n+1}$ isomorphic to $\U n$ which acts transitively on $L$. We claim that there exists $v\in \C^{n+1}$, $v\neq 0$, such that $K = \{g\in \SU{n+1}|\ gv\in \C^* v\}$. Indeed we first note that $K$ acts reducibly on $\C^{n+1}$ because otherwise the center of $K\subset \SU{n+1}$ would act as a multiple of the identity by Schur's Lemma and therefore it would be finite. The claim now follows from the fact that any
irreducible representation of $\SU n$ has dimension at least $n$.

Note that the semisimple part of $H$ is contained in~$K$, and non-trivial
because $n\geq 3$. Therefore $v\in \mathrm{Span}\{e_1,e_2\}$, where $\{e_1,\ldots,e_{n+1}\}$ is the canonical basis of $\C^{n+1}$. Moreover $H$ is not contained in $K$, so  $v\not\in \C e_1$ and $v\not \in \C e_2$. This implies that $H \cap K = \{(z,z,A)\in \mathrm{S}(\U 1\times\U 1\times\U {n-1})\}$.
The isotropy representation of $H$ restricted to $H \cap K$ preserves the subspaces $V_2, V_3$ endowed with the invariant complex structure $J$ and it is of real type on $V_2\oplus V_3$. We can suppose that $J$ is either $J_0$ or $J_1$ and therefore $(V_2,J)$ and $(V_3,J)$ are either equivalent or dual to each other as $H\cap K$-modules. Since the center of $H\cap K$ acts as a scalar multiple of the identity on $V_2\oplus V_3$ and it is of real type, we see that
$(V_2,J)$ and $(V_3,J)$ are dual to each other and thus $J=J_1$. \par

The K\"ahler metric $g$ induces a $J$-Hermitian scalar product on $\bigoplus_{i=1}^3 V_i$
such that this is an orthogonal decomposition. We denote by $g_i$ the
restriction of $g$ on $V_i$ for $i=1$, $2$, $3$.
As $H\cap K$-modules, we can write $V_3 = V_2^*$ and $g_3 = \alpha\cdot g_2^*$, where $\alpha\in\R^+$ and $g_2^*$ is the Hermitian metric induced by $g_2$ on $V_2^*$. We
may assume that the $H\cap K$-invariant real form $\ell$ of $V_2\oplus V_2^*$ is given by
$\ell = \{(v,v^*)|\ v\in V_2\}$, where $v^*\in V_2^*$ is the dual of $v\in V_2$
with respect to~$g_2$. Writing the condition that $\ell$ is Lagrangian
relative to the metric $g$, we immediately see that $\alpha = 1$. This completely determines the metric by the K\"ahler condition (see e.g.~\cite[Theorem~9.4(2)]{W-G}) and it turns out that the projection $\pi:(M,g)\to W$ is a Riemannian submersion when we choose a suitable multiple of the symmetric metric on $W$. This means that the metric $g$ is (up to a multiple) the standard K\"ahler-Einstein metric on the twistor space (see e.g.~\cite[14.80]{Be}). \EPf

\medskip

Now Lemma~\ref{twistor} shows that the standard construction described above
provides the twistor space $Z$ with a tight Lagrangian sphere
$L\cong\U n/\U{n-1}\cong\ S^{2n-1}$.
It is interesting to remark that $L$ coincides with a
\emph{natural lift} in~$Z$
of a totally complex projective space $\C P^{n-1}\subset W$~\cite{ET}.

\smallskip

(b) $n=2$. The flag manifold $M$ is the full flag manifold $\SU 3/T^2$ which admits only one invariant complex structure $J$, up to equivalence. The standard
construction shows that $M$, endowed with a suitable invariant
K\"ahler metric $g$, admits a tight Lagrangian sphere $L\cong S^3$ which is
an orbit of a subgroup $K\cong \U 2$ of $\SU 3$
and a connected component of an antiholomorphic isometry~$\tau$.
Arguments like those in the proof of Lemma~\ref{twistor} show
that $g$ is K\"ahler-Einstein.

\smallskip

We are left with proving that the tight Lagrangian spheres are unique up to the $G$-action. This is clear for the real form $S^n$ in the quadric $Q_n$, so we will focus on the remaining cases. If $L'$ is another tight Lagrangian sphere, we know that $L'$ is homogeneous under the action of a subgroup $K'$ with $(G,K')$ a symmetric pair of rank one. This implies that $K'$ is unique up to conjugation and we can suppose that $(G,K') = (G,K)$, where $(G,K)$ is one of the standard pairs $(\SU{n+1},\U n)$, $(\SP{n+1},\SP 1 \times\SP{n})$ or $(\mathrm{F}_4, \Spin 9)$, respectively.  Therefore it is enough to show that the standard subgroup $K$ has a unique orbit which is a Lagrangian sphere, and this follows
from~\cite[Theorem~1.2]{BG}.


\providecommand{\bysame}{\leavevmode\hbox to3em{\hrulefill}\thinspace}
\providecommand{\MR}{\relax\ifhmode\unskip\space\fi MR }
\providecommand{\MRhref}[2]{%
  \href{http://www.ams.org/mathscinet-getitem?mr=#1}{#2}
}
\providecommand{\href}[2]{#2}

\end{document}